\documentclass{article}
\usepackage{amssymb,amsmath}
\usepackage{graphicx}
\usepackage{ dsfont }

\def\qed{\hfill {\hbox{${\vcenter{\vbox{               %HOLLOW SQUARE
   \hrule height 0.4pt\hbox{\vrule width 0.4pt height 6pt
   \kern5pt\vrule width 0.4pt}\hrule height 0.4pt}}}$}}}

\def\tr{\triangleright}

\newtheorem{theorem}{Theorem}
\newtheorem{definition}{Definition}

\newtheorem{example}{Example}
\newtheorem{remark}[example]{Remark}

{\qed\par\medskip}

\date{}

\title{\Large \textbf{Enhancements of the rack counting 
invariant via $N$-reduced dynamical cocycles}}

\author{Alissa S. Crans\footnote{acrans@lmu.edu} 
\and Sam Nelson\footnote{knots@esotericka.org}
\and Aparna Sarkar\footnote{as072010@mymail.pomona.edu}} 

\begin{document}
\maketitle

\begin{abstract}
We introduce the notion of \textit{$N$-reduced dynamical cocycles}
and use these objects to define enhancements of the rack counting 
invariant for classical and virtual knots and links. We provide 
examples to show that the new invariants are not determined by the 
rack counting invariant, the Jones polynomial or the generalized 
Alexander polynomial.
\end{abstract}



\section{\large \textbf{Introduction}}

\textit{Racks} were introduced in 1992 in \cite{FR} as an algebraic structure
for defining representational and functorial invariants of framed 
oriented knots and links.  A rack generalizes the notion of a  \textit{quandle}, an algebraic 
structure defined in 1980 in \cite{J,M} which defines invariants of unframed 
knots and links. More precisely, the number of quandle homomorphisms from the 
fundamental quandle of a knot or link to a finite quandle $X$ defines a 
computable integer-valued invariant of unframed oriented knots and links known 
as the \textit{quandle counting invariant}.

In \cite{N}, a property of finite racks known as \textit{rack rank} or 
\textit{rack characteristic} was used to define an integer-valued invariant of
unframed oriented knots and links using non-quandle racks, known as the 
\textit{integral rack counting invariant}; for quandles, this invariant
coincides with the quandle counting invariant. An \textit{enhancement} of a 
counting invariant uses a Reidemeister-invariant signature for each 
homomorphism rather than merely counting homomorphisms. In \cite{CJKLS}, the 
first enhancement of the quandle counting invariant was defined using 
\textit{Boltzmann weights} determined by elements of the second cohomology 
of a finite quandle. The resulting \textit{quandle 2-cocycle invariants} of
knots and links have been the subject of much study ever since.

In \cite{HHNYZ} an enhancement of the integral rack counting invariant
was defined using a modification of the \textit{rack module} structure from
\cite{AG}, associating a vector space or module to each homomorphism. In this
paper we further generalize the enhancement from \cite{HHNYZ} using a 
modified version of an algebraic structure first defined in \cite{AG} known 
as a \textit{dynamical cocycle}. In particular, dynamical cocycles satisfying
a condition we call \textit{$N$-reduced} yield an enhancement of the rack
counting invariant.

The paper is organized as follows. In Section \ref{rb} we review the basics of 
racks, the rack counting invariant, and the rack module enhancement. In Section
\ref{dc} we define $N$-reduced dynamical cocycles and the $N$-reduced
dynamical cocycle invariant. In Section \ref{cex} we provide some computations
and examples, and we conclude in Section \ref{q} with some questions for
future study.

\section{\large \textbf{Racks, the counting invariant and the rack module 
enhancement}} \label{rb}

We start by reviewing some basic definitions from \cite{FR,J}.

\begin{definition} 
\textup{A {\it rack } is a set $X$ equipped with a binary operation 
$\tr : X \times X \to X$ satisfying the following two conditions:
\begin{enumerate}
\item[(i)] For each $x\in X$, the map $f_x: X \to X$ defined by 
$f_x(y) = y \rhd x$ is invertible, with inverse $f^{-1}_{x}(y)$ denoted by
$y\rhd^{-1}x$, and
\item[(ii)] For each $x,y,z\in X$, we have 
$(x\rhd y)\rhd z = (x \rhd z) \rhd (y \rhd z).$
\end{enumerate}
A \textit{quandle} is a rack with the added condition:
% that 
\begin{itemize} \item[(iii)] For all $x \in X$, we have $x \rhd x =x$. 
\end{itemize}}
\end{definition}

Note that axiom (ii) is equivalent to the requirement that each map 
$f_x:X\to X$ is a rack homomorphism, i.e.
\[f_z(x\tr y)=(x\tr y)\tr z = (x\tr z)\tr(y\tr z)=f_z(x)\tr f_z(y),\]
so we can alternatively define a rack as a set $X$ with a bijection 
$f_x:X\to X$ for each $x\in X$ such that every $f_x$ is an automorphism
of the structure on $X$ defined by $x\tr y=f_y(x)$.

Standard examples of racks include:
\begin{itemize}
\item \textit{$(t,s)$-racks.} Any module over 
$\ddot\Lambda=\mathbb{Z}[t^{\pm 1},s]/(s^2-(1-t)s)$ is a rack under
\[x\tr y=tx+sy.\]
If $s$ is invertible, then $s^2-(1-t)s=0$ implies $s=1-t$ and we have a
quandle known as an \textit{Alexander quandle}.
\item \textit{Conjugation racks.} Every group $G$ is a rack (indeed, a quandle)
under $n$-fold conjugation for any $n\in \mathbb{Z}$:
\[x\tr y=y^{-n}xy^n.\]
\item \textit{The Fundamental Rack of a framed oriented link.} Let 
$L\subset S^3$ be a link of $c$ components, $n(L)$ a regular neighborhood of 
$L$ with set of framing curves $F=\{F_1,\dots, F_c\}$ giving the framing of 
$L$, $x_0\in S^3\setminus n(L)$ a base point and $FR(L)$ the set of isotopy 
classes of paths from $x_0$ to $F_i$ where the terminal point of the path 
can wander along $F_i$ during the isotopy. For each point $x_1\in F_i$ there 
is a meridian $m(x_1)$ in $n(L)$, unique up to isotopy, linking the $i$th 
component of $L$ once. Then for each path $y:[0,1]\to S^3\setminus n(L)$ 
representing an isotopy class in $FR(L)$, let 
$p(y)=y^{-1}\ast m(y(1))\ast y\in \pi_1(S^3\setminus n(L))$ 
where $\ast$ is path  concatenation reading right-to-left.  Then $FR(L)$ is 
a rack under the operation \[[x]\tr [y]= [x\ast p(y)].\]
Combinatorially, $FR(L)$ can be understood 
as equivalence classes of rack words in a set of generators corresponding 
one-to-one with the set of arcs in a diagram of $L$ under the equivalence r
elation generated by the rack axioms and crossing relations in $L$.
See \cite{FR} for more details.
\end{itemize}

\begin{definition}\textup{
Let $X=\{x_1,\dots,x_n\}$ be a finite set. We can specify a rack structure
on $X$ by a \textit{rack matrix} $M_X$ in which the $(i,j)$th entry is $k$ when
$x_k=x_i\tr x_j$. Rack axiom (i) is equivalent to the condition that every
column of $M_X$ is a permutation; rack axiom (ii) requires checking
each triple for the condition $M_{M_{i,j},k}=M_{M_{i,k},M_{j,k}}$.}
\end{definition}

\begin{example} \label{exM}
\textup{The $(t,s)$-rack structure on $\mathbb{Z}_4=\{x_1=1,x_2=2,x_3=3,x_4=0\}$with $t=1$ and $s=2$ has rack matrix}
\[M_X=\left[\begin{array}{cccc}
3 & 1 & 3 & 1 \\ 
4 & 2 & 4 & 2 \\
1 & 3 & 1 & 3 \\
2 & 4 & 2 & 4 \\
\end{array}\right].\]
\end{example}

\begin{definition}\textup{
Let $X$ be a rack and $L$ an oriented link diagram. An \textit{$X$-labeling} 
or \textit{rack labeling of $L$ by $X$} is an assignment of an element of 
$X$ to each arc in $L$ such that the condition below is satisfied:}
\[\includegraphics{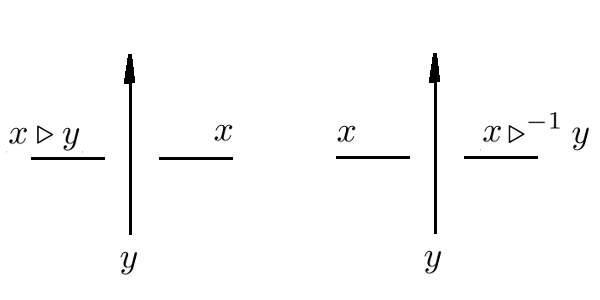}\]
\end{definition}

Indeed, the rack axioms are algebraic distillations of Reidemeister 
moves II and III under this labeling scheme; the quandle condition 
corresponds to the unframed Reidemeister move I, and  the framed Reidemeister 
I moves do not impose any additional conditions. Accordingly, labelings of arcs 
of oriented framed knot or link diagrams by rack elements (respectively, 
quandle elements) as shown above are preserved by oriented framed Reidemeister 
moves (respectively, oriented unframed Reidemeister moves) as illustrated in the figures below.
\[\raisebox{-0.65in}{\includegraphics{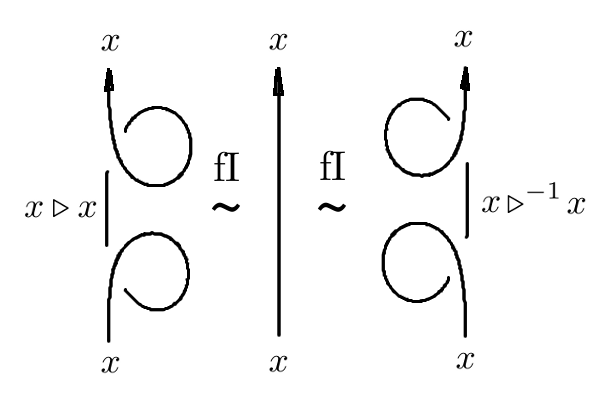}} \quad \left(\mathrm{respectively,}\quad \raisebox{-0.65in}{\includegraphics{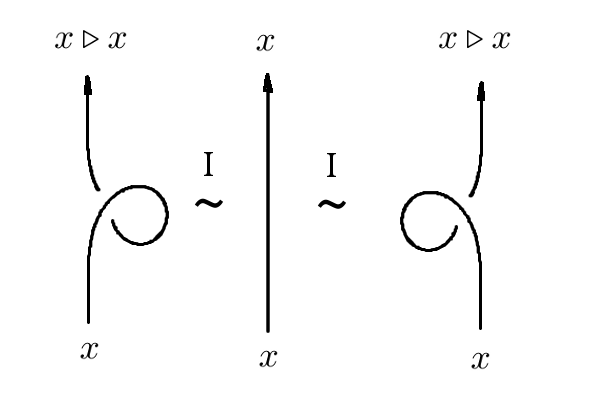}}\right)\]
\[\includegraphics{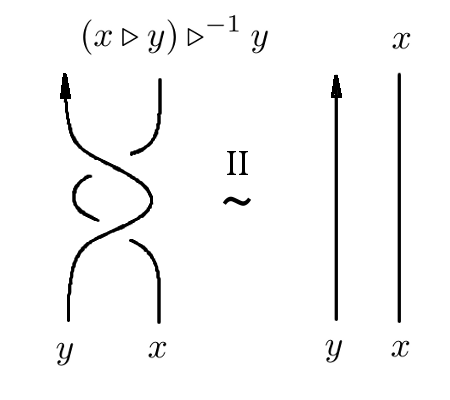} \quad\quad \includegraphics{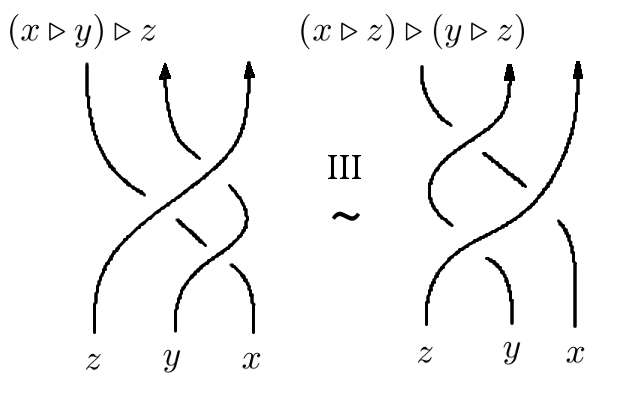}\]

\begin{definition} \textup{Let $X$ be a rack. We call the map $\pi:X \to X$ 
defined by $\pi(x) = x \rhd x$ the \textit{kink map}. The \textit{rack rank} or 
\textit{rack characteristic} of $X$, denoted by $N(X)$, is the order of the 
permutation $\pi$ considered as an element of the symmetric group $S_{|X|}$.
Equivalently, for every element $x\in X$, the \textit{rank} of $x$, denoted by  
$N(x)$, is the smallest positive integer $N$ such that $\pi^N(x) = x.$
Thus, $N(X)$ is the least common multiple of the ranks $N(x)$ for all $x \in X$.
In particular, the kink map of a rack structure on a finite set 
$X=\{x_1,\dots,x_n\}$ given by a rack matrix $M_X$ is 
the permutation in $S_{|X|}$ which sends $k$ to the $(k,k)$ entry of $M_X$.
That is, the image of $\pi$ is given by the entries along the diagonal 
of $M_X$.}
\end{definition}

\begin{example}\textup{
The rack in Example \ref{exM} has kink map satisfying $\pi(1)=3$, $\pi(2)=4$,
$\pi(3)=1$ and $\pi(4)=2$
(or, in cycle notation, $\pi=(13)(24)$) and hence has rack rank $N=2$.
}\end{example}

\begin{remark}\textup{The quandle condition implies that the rank of 
every quandle element is 1, and thus the rack rank of a quandle is 
always 1. Indeed, quandles are simply racks with rack rank $N=1$.}
\end{remark}

Rack rank can be understood geometrically in terms of the Reidemeister type I 
move: if an arc in a knot diagram is labeled with a rack element $x$, going 
through a positive kink changes the label to $\pi(x)$. A natural question is
then: how many kinks must we go though to end up again with $x$? This notion of 
order is the rank of $x$. We can illustrate the concept of rack rank with the 
\textit{$N$-phone cord move} pictured below:
\[\includegraphics{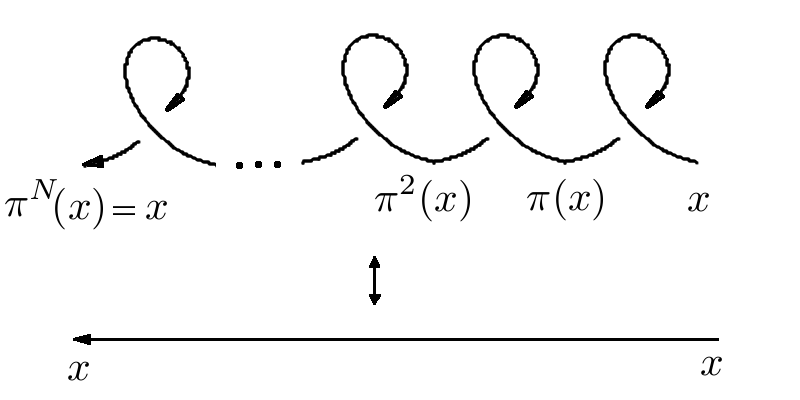} \]
If $N$ is the rank of $X$, then labelings of a link diagram
$L$ by $X$ are preserved by $N$-phone cord moves. In particular, if $X$ is a
rack of rack rank $N$, and $L$ and $L'$ are framed oriented links related by 
framed Reidemeister moves with framings congruent modulo $N$, then
the sets of $X$-labelings of $L$ and $L'$ are in bijective correspondence. 
It follows that the number of homomorphisms is periodic in $N$ on each 
component of a link $L$. 

\begin{definition}\textup{
Let $X$ be a rack with rank $N$ and let $L$ be an oriented link of $c$ 
components. Let $\mathbf{w}\in(\mathbb{Z}_N)^c$ be a framing vector
specifying a framing modulo $N$ for each component of $L$, and let us denote
a diagram of $L$ with framing vector $\mathbf{w}$ by $(L,\mathbf{w})$.  
We thus obtain a set of $N^c$ diagrams of framings of  $L$ mod $N$. For each 
such diagram $(L,\mathbf{w})$, we have a set of $X$-labelings corresponding
to homomorphisms $f:FR(L,\mathbf{w})\to X$.  Summing the numbers of 
$X$-labelings over the set 
$\{(L,\mathbf{w})\ |\ \mathbf{w}\in (\mathbb{Z}_N)^c\}$, we obtain an
invariant of unframed links known as the
\textit{integral rack counting invariant,} which is denoted by:}
\[\Phi_X^{\mathbb{Z}}(L)=\sum_{\mathbf{w}\in(\mathbb{Z}_N)^c} 
|\mathrm{Hom}(FR(L,\mathbf{w}),X)|.\]
\end{definition}

\begin{example}\label{exh}
\textup{Let $X$ be the rack with rack matrix 
$M_X=\left[\begin{array}{cc}
2 & 2  \\
1 & 1
\end{array}\right]$. As a labeling rule, the rack structure of $X$ says
that at a crossing, the understrand switches from 1 to 2 or from 2 to 1
since $1\tr x=2$ and $2\tr x=1$ for $x=1,2$.
The kink map is the transposition $(12)$, so $N=2$. Thus, to compute 
$\Phi_X^{\mathbb{Z}}$ on a link of $c=2$ components, we must count 
$X$-labelings
on the set of $N^c=2^2=4$ diagrams with writhe vectors in $(\mathbb{Z}_N)^c$.
The $(4,2)$-torus link $L4a1$ and the Hopf link $L2a1$ both have
four $X$-labelings as depicted below, so we have 
$\Phi_X^{\mathbb{Z}}(L4a1)=\Phi_X^{\mathbb{Z}}(L2a1)=4$.}
\[\includegraphics{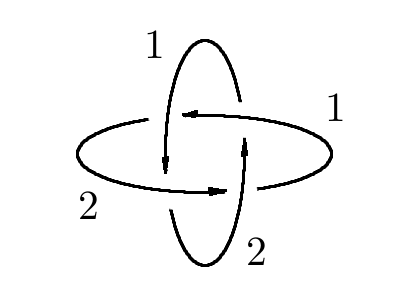}
\includegraphics{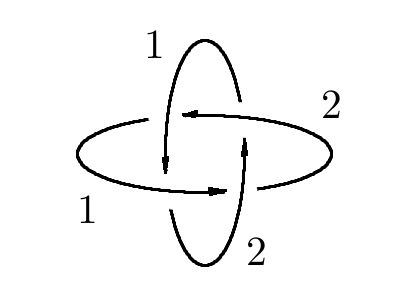}
\includegraphics{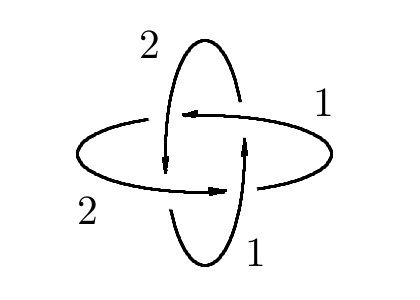}
\includegraphics{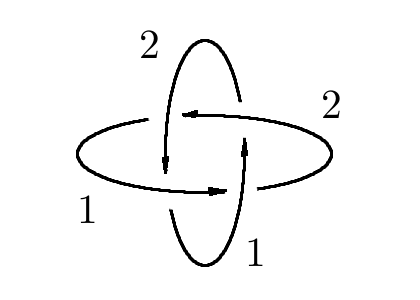}\]
\[
\includegraphics{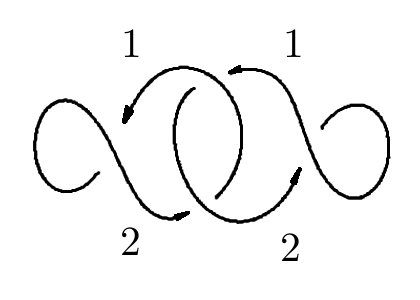}
\includegraphics{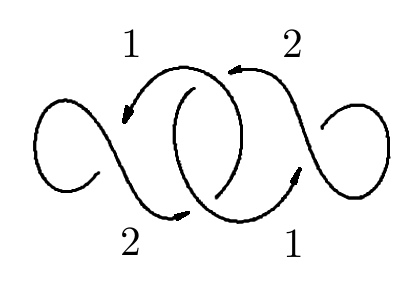}
\includegraphics{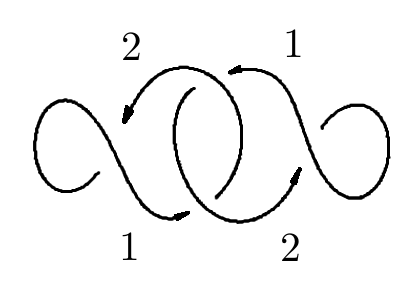}
\includegraphics{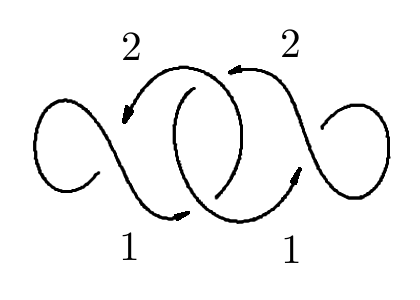}\]
\end{example}

An \textit{enhancement} of $\Phi_{X}^{\mathbb{Z}}(L)$ is a link invariant
defined by associating to each $X$-labeling of $L$ a quantity which is 
unchanged by $X$-labeled framed Reidemeister moves and $N$-phone cord moves. 
Examples include:
\begin{itemize}
\item \textit{Image Enhanced Invariant.} The image of rack homomorphism is
closed under $\tr$ and thus is unchanged by $N$-phone cord moves.  
Hence we have an enhancement:
\[\Phi_X^{\mathrm{Im}}(L)=\sum_{\mathbf{w}\in(\mathbb{Z}_N)^c} 
\left(\sum_{f\in\mathrm{Hom}(FR(L,\mathbf{w}),X)}u^{|\mathrm{Im}(f)|}\right)\]
where $u$ is a formal variable.
\item \textit{Writhe Enhanced Invariant.} Keeping track of which labelings
are contributed by which writhes yields another enhancement:
\[\Phi_X^{\mathrm{W}}(L)=\sum_{\mathbf{w}\in(\mathbb{Z}_N)^c} 
|\mathrm{Hom}(FR(L,\mathbf{w}),X)|q^{\mathbf{w}}\]
where $q^{(w_1,\dots,w_c)}=q_1^{w_1}\dots q_c^{w_c}$ is a product of formal 
variables.
\item \textit{Cocycle Invariants.} A finite rack $X$ has a cohomology theory
analogous to group cohomology. For any 
$f\in\mathrm{Hom}_{\mathbb{Z}}(\mathbb{Z}[X^n],\mathbb{Z})$,  define 
$\delta^n:\mathbb{Z}[X^n]\to \mathbb{Z}[X^{n+1}]$ by
\[\begin{array}{rcl}(\delta^{n} f)(x_1,\dots,x_{n+1}) & = & 
\displaystyle{\sum_{k=2}^{n+1}} 
(-1)^k(f(x_1,\dots,x_{k-1},x_{k+1},\dots, x_{n+1}) \\
& & \quad \quad \quad -f(x_1\tr x_k,\dots,x_{k-1}\tr x_k,x_{k+1},\dots, x_{n+1}))
\end{array}\] and extend linearly. Let $D^n$ be the subgroup of 
$\mathbb{Z}[X^n]$ generated by elements of the form 
\[\sum_{k=1}^N (x_1,\dots, \pi^k(x_j),\pi^{k+1}(x_j),\dots, x_n),\quad 
\quad \quad j=1,\dots, n-1\] 
where $N$ is the rack rank of $X$. Then $(D^n,\delta^n)$ is a subcomplex
of $(\mathbb{Z}[X^n],\delta^n)$; the quotient complex 
$(\mathbb{Z}[X^n]/D^n,\delta^n)$ is the \textit{$N$-reduced rack cochain 
complex} (or the \textit{quandle cochain complex} if $N=1$), with cohomology 
groups denoted by $H^n_{R/ND}(X)$. For every element $\phi\in H^{2}_{R/ND}(X)$ 
(such a $\phi$ is called an\textit{ $N$-reduced 2-cocycle}) we have an 
enhancement
\[\Phi_X^{\phi}(L)=\sum_{\mathbf{w}\in(\mathbb{Z}_N)^c} 
\left(\sum_{f\in\mathrm{Hom}(FR(L,\mathbf{w}),X)}u^{BW(f)}\right)\]
where $BW(f)$, the \textit{Boltzmann weight} of $f$, is the sum over
all crossings in $f$ of $\phi$ evaluated at the arc labelings of each 
crossing.
\end{itemize}
See \cite{CJKLS,EN,N} for further details.

\begin{example}
\textup{In Example \ref{exh}, the links $L2a1$ and $L4a1$ have the same 
number of $X$-labelings over a complete period of framings mod $N$, but
these labelings occur at different framing vectors. In particular, all four 
labelings of $L4a1$ occur with writhe vector $\mathbf{w}=(0,0)$ while all 
four labelings of $L2a1$ occur with writhe vector $\mathbf{x}=(1,1)$.  Thus
the writhe enhanced invariant $\Phi_{X}^W$ distinguishes the links, with
$\Phi_X^W(L4a1)=4\ne4q_1q_2=\Phi_X^W(L2a1)$.}
\end{example}

In \cite{AG} an algebra known as the \textit{rack algebra} $\mathbb{Z}[X]$ 
was associated to each finite rack $X$; in \cite{HHNYZ} a modified form
of the rack algebra was used to define an enhancement of $\Phi_X^{\mathbb{Z}}$.
The idea is to add a secondary labeling to an $X$-labeled link diagram by
putting \textit{beads} on each arc and defining a $(t,s)$-rack style 
operation on the beads at a crossing with $t$ and $s$ values indexed by the
arc labels in $X$ as depicted below:
\[\raisebox{-0.5in}{\includegraphics{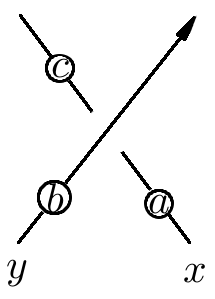}} 
\quad c = t_{x,y}a+s_{x,y}b\]

\begin{definition} 
\textup{Let $X$ be a finite rack with rack rank $N$. The \textit{rack algebra} 
of $X$, denoted by $\mathbb{Z}[X]$, is the quotient of the polynomial 
algebra $\mathbb{Z}[t^{\pm 1}_{x,y}, s_{x,y}]$ generated by noncommuting
variables $t^{\pm 1}_{x,y}$ and $s_{x,y}$ for each
$x,y\in X$ modulo the ideal $I$ generated by the relators
\[ t_{x\tr y,z}t_{x,y} - t_{x\tr z, y \tr z}t_{x,z}, \quad  
t_{x\tr y,z}s_{x,y} - s_{x\tr z, y \tr z}t_{x,z},\quad 
s_{x\tr y,z}- s_{x\tr z, y \tr z}s_{y,z}- t_{x\tr z, y \tr z}s_{x,z}\]
\[ \mathrm{and} \quad \displaystyle{1-\prod^{N-1}_{k=0}\left(t_{\pi^k(x),\pi^k(x)} +s_{\pi^k(x),\pi^k(x)}\right)}\]
for all $x,y,z\in X$.
An \textit{$X$-module} is a representation of $\mathbb{Z}[X]$, that is, an 
abelian group $G$ with automorphisms $t_{x,y}:G\to G$ and endomorphisms
$s_{x,y}:G\to G$ such that the maps defined by the relators of $I$ are zero.} 
\end{definition}

\begin{example}\textup{
Let $R$ be a commutative ring. Then any $R$-module becomes an $X$-module
with a choice of automorphisms and endomorphisms given by multiplication by 
invertible elements $t_{x,y}\in R$ and generic elements 
$s_{x,y}\in R$ such that the ideal $I$ is zero.  We can express such a structure
conveniently with a block matrix $M_R=[\begin{array}{c|c} T & S \end{array}]$ 
where the $(i,j)$ entries of
$T$ and $S$ are $t_{x_i,y_j}$ and $s_{x_i,y_j}$ respectively.}
\end{example}

\begin{example}\textup{
Let $X$ be a rack and let $f\in\mathrm{Hom}(FR(L),X)$ be an $X$-labeled link 
diagram. The \textit{fundamental $\mathbb{Z}[X]$-module} of $f$, denoted by 
$\mathbb{Z}[f]$, is the quotient of the
free $\mathbb{Z}[X]$-module generated by the set of arcs in $f$ modulo the 
ideal generated by the crossing relations.
}\end{example}

In \cite{HHNYZ} an enhancement of $\Phi_X^{\mathbb{Z}}$ was defined using
the number of bead labelings of an $X$-labeled diagram of a framed oriented 
link $L$ as a signature as follows:
\begin{definition}\textup{
Let $X$ be a finite rack and $R$ a commutative ring with an $X$-module
structure. The \textit{rack module enhanced invariant} is given by:
\[\Phi_{X,R}(L)=\sum_{\mathbf{w}\in(\mathbb{Z}_N)^c}\left(\sum_{f\in\mathrm{Hom}(FR(L,\mathbf{w}),X)}u^{|\mathrm{Hom}(\mathbb{Z}[f],R)|}\right).\]
}\end{definition}

\begin{example}
\textup{Let $X$ be the rack from Example \ref{exh} and let $R=\mathbb{Z}_3$. 
The matrix
\[M_R=[T|S] = \left[ {\begin{array}{cc|cc}
 1 & 1 & 1 & 2 \\
 1 & 1 & 2 & 1  \\
 \end{array} } \right]\] defines an $X$-module structure on $R$. 
To compute $\Phi_{X,R}$ for the Hopf link $L2a1$, we must compute 
$|\mathrm{Hom}(\mathbb{Z}[f],R)|$ for each valid $X$-labeling of 
$L2a1$. For instance, the following $X$-labeled diagram has fundamental
$\mathbb{Z}[X]$-module with listed presentation matrix:
\[
\raisebox{-0.5in}{\includegraphics{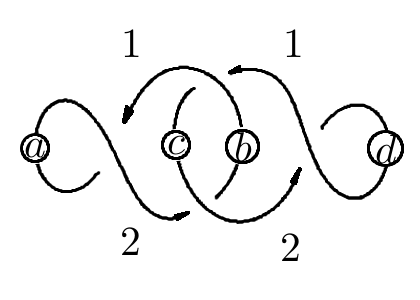}}\qquad
M_{\mathbb{Z}[f]}=\left[\begin{array}{cccc}
t_{2,2}+s_{2,2} & -1 & 0 & 0 \\
0 & s_{1,1} & -1 & t_{1,1} \\
t_{2,2} & -1 & s_{2,2} & 0 \\
0 & 0 & -1 & t_{1,1}+s_{1,1} \\ 
\end{array}\right]\]
Replacing each $t_{x,y}$ and $s_{x,y}$ with its value from $M_R$ and 
row-reducing over $\mathbb{Z}_3$, we have
\[\left[\begin{array}{cccc}
2 & 2 & 0 & 0 \\
0 & 1 & 2 & 1 \\
1 & 2 & 1 & 0 \\
0 & 0 & 2 & 2 \\ 
\end{array}\right]\rightarrow
\left[\begin{array}{cccc}
1 & 0 & 0 & 1 \\
0 & 1 & 0 & 2 \\
0 & 0 & 1 & 1 \\
0 & 0 & 0 & 0 \\ 
\end{array}\right],
\]
so the solution space (i.e., the set of bead labelings) is the set
$\{(0,0,0,0),(2,1,2,1),(1,2,1,2)\}$
and this $X$-labeling contributes $u^3$ to $\Phi_{X,R}(L2a1)$. Repeating for
the other labelings, we have $\Phi_{X,R}(L2a1)=4u^3$.}\end{example}

\section{\large \textbf{Dynamical cocycles and enhancements of 
the counting invariant}}\label{dc}

In this section we generalize the rack module idea to remove the
restrictions of the abelian group structure, keeping only those conditions
required by the Reidemeister moves. The result is a rack structure on 
the product $X\times S$ defined via a map 
$\alpha:X\times X\to \mathrm{Maps}(S\times S,S)$ known as a \textit{dynamical 
cocycle}. Dynamical cocycles were defined in \cite{AG} and used to 
construct extension racks; we will use dynamical cocycles satisfying
an extra condition, which we call \textit{$N$-reduced dynamical cocycles}, to
define an enhancement of the rack counting invariant $\Phi_X^{\mathbb{Z}}$.  

\begin{definition} \textup{Let $X$ be a finite rack of rack rank $N$ and $S$ be 
a finite set.  The elements of $S$ will be called \textit{beads}. A 
map $\alpha: X \times X \to \mathrm{Maps}(S\times S,S)$ may be understood 
as a collection of binary operations $\cdot_{x,y}:S\times S\to S$ indexed by 
pairs of elements of $X$ where where we write 
$a\cdot_{x,y} b=\alpha(x,y)(a,b)$.  Such a map $\alpha$ is a
\textit{dynamical cocycle} on $S$ if the maps satisfy: 
\begin{itemize}
\item[(i)] For all $x,y\in X$ and $b\in S$, the map $f_b^{x,y}:S\to S$ defined by
$f_b^{x,y}(a)=a\cdot_{x,y}b$ is a bijection, and
\item[(ii)] For all $x,y,z\in X$ and $a,b,c\in S$, we have
\[(a\cdot_{x,y}b)\cdot_{x\tr y,z} c=(a\cdot_{x,z}c)\cdot_{x\tr z,y\tr z}(b\cdot_{y,z} c).\]
\end{itemize}}
\end{definition}

\begin{definition}
\textup{Let $X$ be a rack of rack rank $N$ and 
$\alpha:X\times X\to \mathrm{Maps}(S\times S,S)$ a dynamical cocycle.
Define $\rho_x:S\to S$ by $\rho_x(a)=a\cdot_{x,x} a$. Then if the 
diagram \[\includegraphics{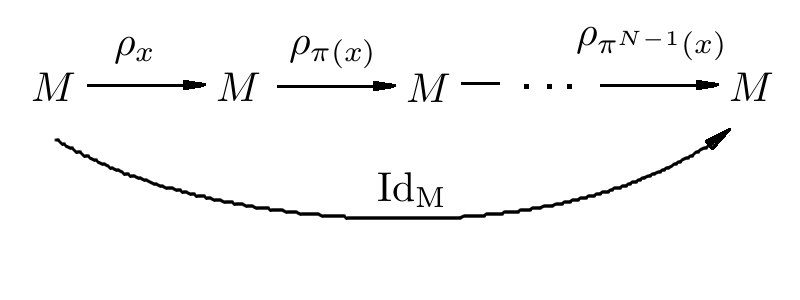}\]
commutes for every $x\in X$ and $a\in S$, we say the cocycle $\alpha$ is 
\textit{$N$-reduced.}}
\end{definition}

The definition of a dynamical cocycle is chosen so that bead labelings 
of an $X$-labeled diagram are preserved under $X$-labeled framed oriented
Reidemeister moves as shown below: 
\[\begin{array}{l} \quad \quad \includegraphics{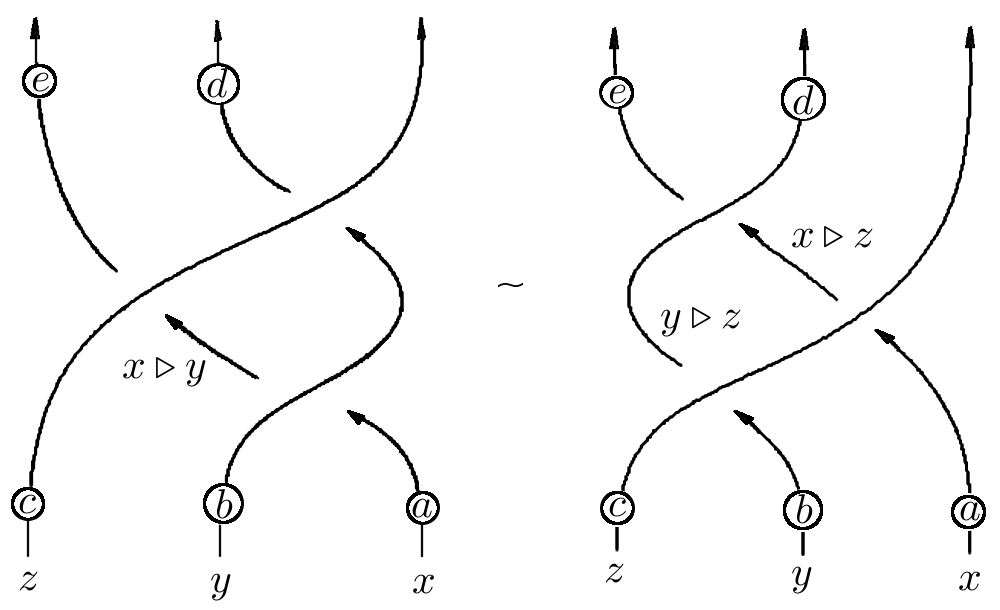} \\
\begin{array}{rclcrcl}
d & = & b \cdot_{y,z} c & \hskip 0.5in  & d & = & b \cdot_{y,z} c\\
e & = & (a\cdot_{x,y}b)\cdot_{x\tr y,z} c & & e & = 
& (a\cdot_{x,z}c)\cdot_{x\tr z,y\tr z}(b\cdot_{y,z} c)\\
\end{array}
\end{array}\]

The Reidemeister II and framed type I moves require the operations
$\cdot_{x,y}:S\times S\to S$ to be right-invertible; the $N$-reduced 
condition is required by the $N$-phone cord move:
\[\raisebox{-0.5in}{\includegraphics{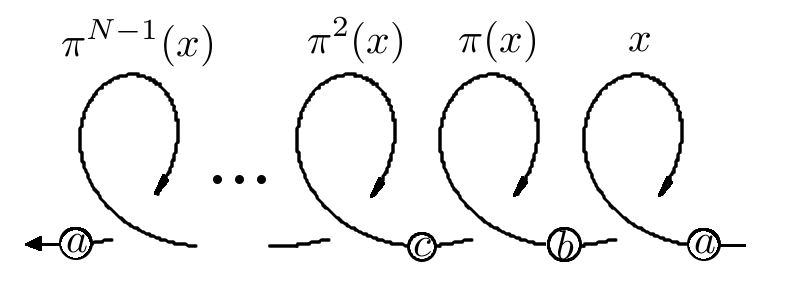}} \quad 
\begin{array}{rcl}
b & = & \rho_{x}(a) \\ 
c & = & \rho_{\pi(x)}(b)=\rho_{\pi(x)}(\rho_x(a)) \\
& \vdots &  \\
a & =& \rho_{\pi^N(x)}(\rho_{\pi^{N-1}(x)}(\dots(\rho_x(a))\dots) \\
\end{array}
\]

\begin{example}
\textup{Let $X$ be a finite rack and $M$ an $X$-module as defined in Section
\ref{rb}. Then the operations
\[a\cdot_{x,y} b=t_{x,y}a+s_{x,y}b\]
define an $N$-reduced dynamical cocycle on $M$.
}\end{example}

More generally, if $X$ is a finite rack of cardinality $n$, we can describe 
a dynamical cocycle on a finite set $S=\{b_1,\dots, b_k\}$ with an 
$(nk)\times (nk)$ block matrix, $M_{x,y},$ encoding the operations tables 
for $\cdot_{x,y}$ 
\[ M_{x,y} = \left[\begin{array}{c|c|c|c}
M_{1,1} & M_{1,2} & \dots & M_{1,n} \\ \hline
M_{2,1} & M_{2,2} & \dots & M_{2,n} \\ \hline
\vdots & \vdots & \ddots & \vdots \\ \hline
M_{n,1} & M_{n,2} & \dots & M_{n,n} \\
\end{array}\right]\]
where the $(i,j)$th entry of $M_{x,y}$ is $l$ 
%such that
when $b_i\cdot_{x,y} b_j=b_l$.

\begin{definition}
\textup{Let $X$ be a finite rack and $\alpha$ an $N$-reduced dynamical cocycle 
on a set $S$. For an $X$-labeled link diagram $f$, let $\mathcal{L}(f)$ be the 
set of $S$-labelings of $f$. Then we define the 
\textit{$N$-reduced dynamical cocycle enhanced invariant} or
\textit{$\alpha$-enhanced invariant} $\Phi_{X,\alpha}(L)$ by:}
\[\Phi_{X,\alpha}(L)=\sum_{\mathbf{w}\in W} 
\left(\sum_{f\in\mathrm{Hom}(FR(L,\mathbf{w}))} u^{|\mathcal{L}(f)|}\right).\]
\end{definition}

By construction, we have
\begin{theorem}
Let $X$ be a finite rack and $\alpha$ an $N$-reduced dynamical cocycle on a 
set $S$. If $L$ and $L'$ are ambient isotopic links, then 
$\Phi_{X,\alpha}(L)=\Phi_{X,\alpha}(L')$.
\end{theorem}

\begin{remark}
\textup{The $\alpha$-enhanced invariant is well-defined for virtual knots
by the usual convention of ignoring virtual crossings.}
\end{remark}

\section{\large \textbf{Computations and Examples}}\label{cex}

In this section we present example computations of the $N$-reduced
dynamical cocycle enhanced invariant.

\begin{example}\label{exJ}
\textup{Let $X$ be the rack with rack matrix 
$M_X=\left[\begin{array}{cc}
2 & 2 \\
1 & 1
\end{array}\right]$
and let $\alpha$ be the dynamical cocycle on $S=\{1,2,3\}$ given by the 
block matrix
\[M_{\alpha}=\left[\begin{array}{ccc|ccc}
3 & 1 & 2 & 2 & 1 & 3 \\ 
1 & 2 & 3 & 3 & 2 & 1 \\
2 & 3 & 1 & 1 & 3 & 2 \\ \hline
2 & 1 & 3 & 3 & 1 & 2 \\
3 & 2 & 1 & 1 & 2 & 3 \\
1 & 3 & 2 & 2 & 3 & 1
\end{array}\right]\]
The virtual knots $3.7$ and the unknot both have Jones polynomial 
$1$ and integral rack counting invariant $\Phi_X^{\mathbb{Z}}=2$. Let us compare 
$\Phi_{X,\alpha}(3.7)$ with $\Phi_{X,\alpha}(\mathrm{Unknot})$.
%but are distinguished by $\Phi_{X,\alpha}$ with
%$\Phi_{X,\alpha}(3.7)=2u^9\ne 2u^3=\Phi_{X,\alpha}(\mathrm{Unknot})$.
Since $X$ has rank $N=2$, we need to consider diagrams with writhes mod 2.
The odd writhe diagrams have no valid $X$-labelings, and there are two valid 
$X$-labelings of the even writhe diagrams. We collect the valid bead labelings
in the tables below.}
\[\raisebox{-0,5in}{\includegraphics{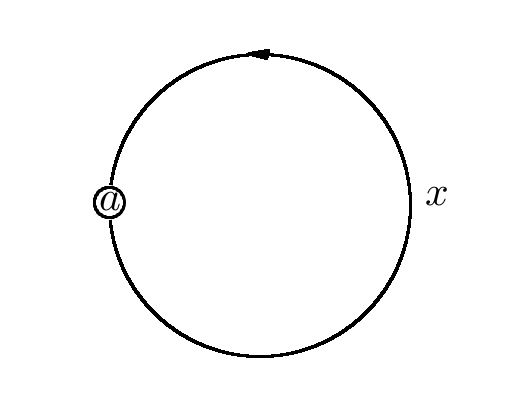}} \quad
\begin{array}{c|c||c|c} 
x & a & x & a \\ \hline
1 & 1 & 2 & 1 \\
1 & 2 & 2 & 2 \\
1 & 3 & 2 & 3 \\
\end{array}\]
\[\raisebox{-0.5in}{\includegraphics{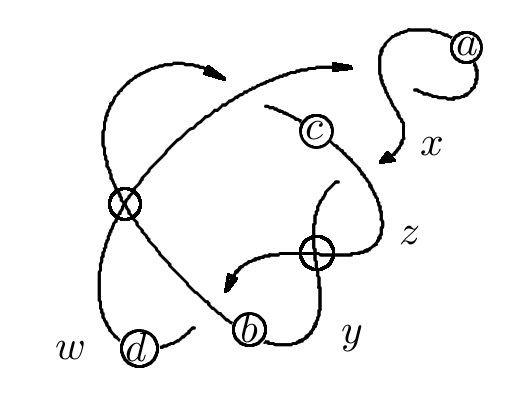} } \quad
\begin{array}{cccc|cccc||cccc|cccc} 
x & y & z & w & a & b & c & d & x & y & z & w & a & b & c & d \\\hline
1 & 2 & 1 & 2 & 1 & 1 & 2 & 3 & 2 & 1 & 2 & 1 & 1 & 1 & 2 & 3 \\
1 & 2 & 1 & 2 & 1 & 2 & 3 & 3 & 2 & 1 & 2 & 1 & 1 & 2 & 3 & 3 \\
1 & 2 & 1 & 2 & 1 & 3 & 1 & 3 & 2 & 1 & 2 & 1 & 1 & 3 & 1 & 3 \\
1 & 2 & 1 & 2 & 2 & 1 & 1 & 2 & 2 & 1 & 2 & 1 & 2 & 1 & 1 & 2 \\
1 & 2 & 1 & 2 & 2 & 2 & 2 & 2 & 2 & 1 & 2 & 1 & 2 & 2 & 2 & 2 \\
1 & 2 & 1 & 2 & 2 & 3 & 3 & 2 & 2 & 1 & 2 & 1 & 2 & 3 & 3 & 2 \\
1 & 2 & 1 & 2 & 3 & 1 & 3 & 1 & 2 & 1 & 2 & 1 & 3 & 1 & 3 & 1 \\
1 & 2 & 1 & 2 & 3 & 2 & 1 & 1 & 2 & 1 & 2 & 1 & 3 & 2 & 1 & 1 \\
1 & 2 & 1 & 2 & 3 & 3 & 2 & 1 & 2 & 1 & 2 & 1 & 3 & 3 & 2 & 1 \\
\end{array}\]
\textup{Hence, we have $\Phi_{X,\alpha}(3.7)=2u^9\ne 2u^3=
\Phi_{X,\alpha}(\mathrm{Unknot})$ and
$\Phi_{X,\alpha}$ is not determined by the Jones polynomial or the integral 
rack counting invariant $\Phi_X^{\mathbb{Z}}$.}
\end{example}
\begin{example}
\textup{Similarly,
the virtual knots $3.7$ and $4.85$ both have generalized
Alexander polynomial 
\[\Delta=(t^2-1)(s^2-1)(st-1)\]
but are distinguished by $\Phi_{X,\alpha}$ with
$\Phi_{X,\alpha}(3.7)=2u^9\ne 2u^3=\Phi_{X,\alpha}(4.85)$ for the rack
$X$ and dynamical cocycle $\alpha$ from Example \ref{exJ}.}
\[\raisebox{-0.5in}{\includegraphics{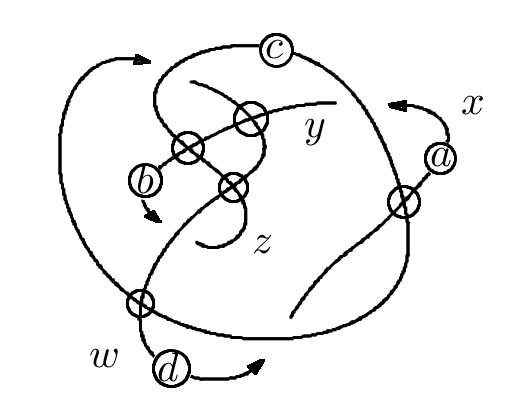} } \quad
\begin{array}{cccc|cccc||cccc|cccc} 
x & y & z & w & a & b & c & d & x & y & z & w & a & b & c & d \\\hline
1 & 2 & 1 & 2 & 1 & 3 & 1 & 3 & 2 & 1 & 2 & 1 & 1 & 3 & 1 & 3 \\
1 & 2 & 1 & 2 & 2 & 2 & 2 & 2 & 2 & 1 & 2 & 1 & 2 & 2 & 2 & 2 \\
1 & 2 & 1 & 2 & 3 & 1 & 3 & 1 & 2 & 1 & 2 & 1 & 3 & 1 & 3 & 1 \\
\end{array}\]
\textup{Hence, $\Phi_{X,\alpha}$ is not determined by the generalized 
Alexander polynomial.}
\end{example}

\begin{example}
\textup{We randomly selected a small dynamical cocycle $\alpha$ on the set
$S=\{1,2,3\}$ for the dihedral quandle $X$ with matrices below.} 
\[M_X=\left[\begin{array}{ccc}
1 & 3 & 2 \\
3 & 2 & 1 \\
2 & 1 & 3 \\
\end{array}\right]
\qquad M_{\alpha}=
\left[\begin{array}{ccc|ccc|ccc}
1 & 3 & 2 & 3 & 2 & 1 & 1 & 3 & 2 \\
3 & 2 & 1 & 2 & 1 & 3 & 3 & 2 & 1 \\
2 & 1 & 3 & 1 & 3 & 2 & 2 & 1 & 3 \\ \hline
3 & 2 & 1 & 1 & 3 & 2 & 2 & 1 & 3 \\
2 & 1 & 3 & 3 & 2 & 1 & 1 & 3 & 2 \\
1 & 3 & 2 & 2 & 1 & 3 & 3 & 2 & 1 \\ \hline
1 & 3 & 2 & 2 & 1 & 3 & 1 & 3 & 2 \\
3 & 2 & 1 & 1 & 3 & 2 & 3 & 2 & 1 \\
2 & 1 & 3 & 3 & 2 & 1 & 2 & 1 & 3 \\
\end{array}\right]
\]
\textup{We then 
computed $\Phi_{X,\alpha}$ for the list of prime classical knots with up to
eight crossings and prime classical links with up to seven crossings
as listed at the \texttt{knot atlas} \cite{KA}. The results are 
collected below. In particular, note that the invariant values 
$6+3u^9\ne 9u^9$ both specailize to the same rack counting invariant
value $\Phi_{X}^{\mathbb{Z}}=9$, and we see that $\Phi_{X,\alpha}$ is
not determined by $\Phi_{X}^{\mathbb{Z}}$.}
\[\begin{array}{r|l}
\Phi_{X,\alpha}(L) & L \\ \hline
3u^3 & \mathrm{Unknot}, 4_1, 5_1, 5_2, 6_2, 6_3, 7_1, 7_2, 7_3, 7_5, 7_6, 8_1,
8_2, 8_3, 8_4, 8_6, 8_7, 8_8, 8_9, 8_{12}, 8_{13}, 8_{14}, 8_{16}, 8_{17}, \\
 & L2a1,
L4a1, L5a1, L6a2, L6a4, L6n1, L7a2, L7a3, L7a4, L7a6, L7a7, L7n1, L7n2 \\
6+3u^9 & 3_1, 7_4, 7_7, 8_5, 8_{15}, 8_{19}, 8_{21}, L6a1, L6a3, L6a5, L7a1 \\
9u^9 & 6_1, 8_{10}, 8_{11}, 8_{20}, L7a5 \\
24+3u^{27} & 8_{18} \\
\end{array}\]

\textup{Our python results indicate that of the 116 prime virtual knots with 
up to 4 classical crossings listed at the knot atlas, $\Phi_{X,\alpha}$ 
for this $\alpha$ is $6+3u^9$ for the virtual knots $3.6, 3.7, 4.61, 4.61, 
4.63, 4.64, 4.65,$ $4.66, 4.67, 4.68 $ and $4.98$, $\Phi_{X,S}^{\alpha}=9u^9$ 
for $4.99$, and $\Phi_{X,\alpha}=3u^3$ for the other virtual knots in 
the list.}
\end{example}

Our \texttt{python} code for computing $N$-reduced dynamical cocycles
and their link invariants is available at \texttt{www.esotericka.org}.

\section{\large \textbf{Questions for future research}}\label{q}

In this section we collect a few questions for future research. 

For a given pair of knots or links, how can we choose $X$ and $\alpha$
to maximize the liklihood of $\Phi_{X,\alpha}$ distinguishing the knots
or links in question? Is there an algorithm, perhaps starting with 
presentations of the fundamental racks of the knots, to construct a rack
$X$ and dynamical cocycle $\alpha$ such that $\Phi_{X,\alpha}$ always
distinguishes inequivalent knots?

A natural direction of generalization is to look at knotted surfaces in
$\mathbb{R}^4$, which have an integral quandle counting invariant which
should be susceptible to enhancement by beads. What analog of the dynamical
cocycle condition arises from the Roseman moves with beads on each sheet?

\bigskip

\noindent\textsc{Department of Mathematics\\ 
Loyola Marymount University \\
One LMU Drive, Suite 2700\\ 
Los Angeles, CA 90045}

\medskip

\noindent\textsc{Department of Mathematics \\ 
Claremont McKenna College \\
850 Columbia Ave. \\
Claremont, CA 91711}

\medskip

\noindent\textsc{Department of Mathematics \\ 
Pomona College \\
610 North CollegeAve \\
Claremont, CA 91711}


\begin{thebibliography}{0}

\bibitem{AG}{N. Andruskiewitsch and M. Gra\~{n}a.
From racks to pointed Hopf algebras. 
\textit{Adv. Math.} \textbf{178} (2003) 177-243.}

\bibitem{CJKLS}{J. S. Carter, D. Jelsovsky, S. Kamada, L. Langford and 
M. Saito. Quandle cohomology and state-sum invariants of knotted curves 
and surfaces.  \textit{Trans. Am. Math. Soc.} \textbf{355} (2003) 3947-3989.}

\bibitem{KA}{D. Bar-Natan (Ed.). The Knot Atlas.
\texttt{http://katlas.math.toronto.edu/wiki/Main\underline{\ }Page}}

\bibitem{CN}{J. Ceniceros and S. Nelson. $(t,s)$-racks and their link 
invariants. arXiv:1011.5455, to appear in \textit{Int'l. J. Math.}}

\bibitem{EN}{M. Elhamdadi and S. Nelson. $N$-degeneracy in rack homology.
arXiv:  to appear in \textit{Hiroshima Math J.}}

\bibitem{FR}{R. Fenn and C. Rourke.
 Racks and links in codimension two.
 \textit{J. Knot Theory Ramifications}  \textbf{1}  (1992) 343-406.}

\bibitem{HHNYZ}{A. Haas, G. Heckel, S. Nelson, J. Yuen, Q. Zhang. 
Rack Module Enhancements of Counting Invariants. arXiv:1008.0114, to appear 
in \textit{Osaka J. Math.}}

\bibitem{J}{D. Joyce.
 A classifying invariant of knots, the knot quandle.
 \textit{J. Pure Appl. Algebra}  \textbf{23}  (1982)  37-65.}

\bibitem{M}{S. V. Matveev.
Distributive groupoids in knot theory.
\textit{Math. USSR, Sb.} \textbf{47} (1984) 73-83.}

\bibitem{N}{S. Nelson. Link invariants from finite racks. arXiv:0808.0029;
to appear in Fund. Math. }



\end{thebibliography}
\end{document}